\documentstyle[epsf]{article}

\makeatletter
  
  \@addtoreset{equation}{section}
\makeatother

\newtheorem{theorem}{Theorem}[section]
\newtheorem{proposition}{Proposition}[section]
\newtheorem{lemma}{Lemma}[section]

\newtheorem{example}{Example}[section]

\newfont{\germ}{eufm10}
\newfont{\germlarge}{eufm10}
\newfont{\slsmall}{cmsl8}

\def\B{{\cal B}}
\def\bbar{\overline{b}}
\def\cd{\cdots}
\def\ch{\mbox{\sl ch}\,}
\def\et#1{\tilde{e}_{#1}}
\def\ft#1{\tilde{f}_{#1}}
\def\geh{\goth{g}}
\def\goth#1{\mbox{\germ #1}}
\def\L{{\cal L}}
\def\La{\Lambda}
\def\la{\lambda}
\def\ot{\otimes}
\def\P{{\cal P}}
\def\pbar{\overline{p}}
\def\Pcl{P_{cl}}
\def\Pcll{(P_{cl}^+)_l}

\def\Q{{\bf Q}}

\def\sln{\goth{sl}_{\,n}}
\def\slnh{\widehat{\goth{sl}}_{\,n}}
\def\Uq{U_q(\geh)}
\def\Uqp{U'_q(\geh)}
\def\veps{\varepsilon}
\def\vphi{\varphi}
\def\wt{\mbox{\sl wt}\,}
\def\Z{{\bf Z}}
\def\Zn{\Z_{\ge0}}

\begin{document}

\title{ Paths, Demazure Crystals \\
        and Symmetric Functions }

\author{
Atsuo Kuniba\thanks{
Institute of Physics, University of Tokyo, Komaba, Tokyo 153, Japan}, 
Kailash C. Misra\thanks{
Department of Mathematics, North Carolina State University, 
Raleigh, NC 27695-8205, USA}, 
Masato Okado\thanks{
Department of Mathematical Science, Faculty of Engineering Science,
Osaka University, Toyonaka, Osaka 560, Japan},\\
Taichiro Takagi\thanks{
Department of Mathematics and Physics, National Defense Academy,
Yokosuka 239, Japan}
and Jun Uchiyama\thanks{
Department of Physics, Rikkyo University, Nishi-Ikebukuro, Tokyo 171, Japan}
} 

\date{}
\maketitle

\begin{abstract}
\noindent
We review the path realization of Demazure crystals and 
discuss Demazure characters in the light of
symmetric functions.
\end{abstract}

\section{Introduction}

Let $U_q(\geh)$ be a quantum affine Lie algebra.
Representation of $U_q(\geh)$ at $q=0$ is well described
by the crystal base theory \cite{Ka1}, \cite{KMN1}, \cite{KMN2}.
For example, consider the irreducible highest weight
$\Uq$-module $V(\lambda)$ for any dominant integral
weight $\lambda$ of level $l$.
At $q=0$ its crystal $\B(\la)$ 
admits a parametrization in terms of {\em paths}.
The latter is the combinatorial object that arose 
in the studies of solvable lattice models 
\cite{DJKMO1}, \cite{DJKMO2} by Baxter's corner transfer
matrix method \cite{Bax}.
Given a perfect crystal $B$ of level $l$, 
a path is an element of the semi-infinite 
tensor product $\cdots B \otimes B$.
It must obey some boundary condition 
on the left tail, which is uniquely specified 
from $\lambda$ and $B$.
Letting $\P(\la, B)$ denote the set of such paths, one has 
an isomorphism of crystals 
$\psi: \B(\la) \buildrel \sim \over \rightarrow \P(\la, B)$.
These features \cite{KMN1}, \cite{KMN2} will be summarized 
in section 2.

In \cite{Ka2}, Kashiwara showed that 
for each Weyl group element $w$ there exists a finite 
subset $\B_w(\la) \subset \B(\la)$ that corresponds to the 
crystal of the Demazure module $V_w(\la) \subset V(\la)$.
Then a natural question arises; What kind of paths are
contained in the image $\psi(\B_w(\la))$ ?
This was answered in \cite{KMOU} for a class of $w$ obeying certain 
conditions.
The result is given for each value of the `mixing index' 
$\kappa \in {\bf Z}_{\ge 1}$ that reflects some property of 
$\lambda$ and $B$.
For simplicity we shall exclusively consider  
$\kappa = 1$ case in this paper 
and refer to \cite{KMOU} for $\kappa$ general case.
Then roughly speaking there occur only those paths corresponding
to a tensor product of some 
subset $B^{(j)}_a \subset B$ and finitely many $B$'s.
The precise description will be given in section 3,
which constitutes the first main contents
in this report.
The result may also be viewed as a combinatorial 
explanation of the tensor product structure 
in the Demazure modules observed in \cite{S}.
In section 4 we shall give an example from 
$\slnh$.
Results on the other classical affine Lie algebras
are available in \cite{KMOTU1}.

Section 5 is devoted to our another main topic in this paper, 
namely, the characters of $V_w(\la)$ 
in the light of symmetric functions.
By Theorem \ref{th:rectangular} the Demazure characters 
provide a $q$-analogue of the products of classical characters,
that is, Schur functions. We then relate them with the Kostka-Foulkes
polynomial and Milne polynomial. Results on some Demazure
characters of the 
other classical affine Lie algebras are available in \cite{KMOTU2}.

A.K. appreciates the great hospitality of
the organizers of the Nankai-CRM meeting
{\em Extended and Quantum Algebras and their Applications to Physics}
during August 19-24, 1996, where a part of this work was presented.
He also thanks T. Nakanishi for hospitality and support at
University of North Carolina.
M.O. thanks O. Foda for hospitality and support at The University 
of Melbourne. He also thanks P. Littelmann, B. Leclerc, J.-Y.
Thibon and T.A. Welsh for discussions.
A.K. and M.O. thank Anatol N. Kirillov for stimulating 
discussions.
K.C.M. is supported in part by NSA/MSP Grant No. 96-1-0013.

\section{Perfect Crystals and Demazure Modules}

First we fix the notations following \cite{KMN1}.
$\Uq$ is the quantized universal enveloping algebra of 
an affine Lie algebra $\geh$. 
Let 
$\{ \alpha_i \}_{i \in I}, \ 
\{ h_i \}_{i \in I} \mbox{ and }  
\{ \Lambda_i \}_{i \in I}$ denote the set of 
simple roots, coroots and fundamental weights.
$P$ is the weight lattice and $P_+=\{\la\in P\mid \langle\la,h_i\rangle\ge0
\mbox{ for any }i\}$. 
$V(\la)$ is the irreducible highest weight module
of highest weight $\la\in P_+$ and 
$(\L(\la),\B(\la))$ is its crystal base
(which was originally denoted by $(L(\la), B(\la))$ in \cite{KMN1}).
For the notation of a finite-dimensional representation of $\Uqp$, we follow 
section 3 in \cite{KMN1}. For instance, $\Pcl$ is the classical weight
lattice, $\Uqp$ is the subalgebra of $\Uq$ generated by $e_i,f_i,q^h$
($h\in(\Pcl)^*$) and $\mbox{Mod}^f(\geh,\Pcl)$ is the category of 
finite-dimensional $\Uqp$-modules which have the weight decompositions.
We set $\Pcl^+=\{\la\in\Pcl\mid\langle\la,h_i\rangle\ge0\mbox{ for any }i\}
\simeq\sum\Zn\La_i$ and 
$(\Pcl^+)_l=\{\la\in\Pcl^+\mid\langle\la,c\rangle=l\}$, 
where $c$ is the canonical central element. Assume $V$ in $\mbox{Mod}^f
(\geh,\Pcl)$ has a crystal base $(L,B)$. For an element $b$ of $B$, we
set $\veps_i(b)=\max\{n\ge0\mid\et{i}^n b\in B\}$, $\veps(b)=\sum_i
\veps_i(b)\La_i$ and $\vphi_i(b)=\max\{n\ge0\mid\ft{i}^n b\in B\}$,
$\vphi(b)=\sum_i\vphi_i(b)\La_i$.

Let $B$ be a perfect crystal of level $l$.
We refer Definition 4.6.1 in \cite{KMN1} for its definition.
For $\la\in(\Pcl^+)_l$, let 
$b(\la)\in B$ be the element such that $\vphi(b(\la))=\la$. 
{}From the definition
of perfect crystal, such a $b(\la)$ exists and is unique. 
Let $\sigma$ be the
automorphism of $(\Pcl^+)_l$ given by $\sigma\la=\veps(b(\la))$. We set 
$\bbar_k=b(\sigma^{k-1}\la)$ and $\la_k=\sigma^k\la$. Then perfectness
assures that we have the isomorphism of crystals
\[
\B(\la_{k-1})\simeq \B(\la_k)\ot B.
\]
Iterating this isomorphism, we have 
\[
\psi_k:\B(\la)\simeq \B(\la_k)\ot B^{\ot k}.
\]
Defining the set of paths $\P(\la,B)$ by
\[
\P(\la,B)=\{p=\cdots\ot p(2)\ot p(1)\mid p(j)\in B,p(k)=
\bbar_k\mbox{ for }k\gg0\},
\]
we see that there is an 
isomorphism of crystals
$\psi: \B(\la) \buildrel \sim \over \rightarrow \P(\la,B)$. 
In particular, the image 
$\psi(u_\la)$ of the 
highest weight vector $u_\la \in \B(\la)$ is given by 
$\pbar = \cd\ot\bbar_k\ot\cd\ot\bbar_2\ot
\bbar_1$.
We call $\pbar$ the {\em ground-state} path.
One can explicitly describe the weights and the actions of 
$\et{i}$ and $\ft{i}$ on $\P(\la,B)$ by 
the energy function and the signature rule.
See sections 1.3 and 1.4 of \cite{KMOU}.

Now we proceed to Demazure crystals.
Let $\{r_i\}_{i\in I}$ be the set of simple reflections, and let
$W$ be the Weyl group. For $w\in W$, $l(w)$ denotes the length of 
$w$, and $\prec$ denotes the Bruhat order on $W$. Let $U_q^+(\geh)$
be the subalgebra of $\Uq$ generated by $e_i$'s. 
For $\la \in (\Pcl^+)_l$ we consider the
irreducible highest weight $\Uq$-module $V(\la)$ as before.
Let $V_w(\la)$ denote the $U_q^+(\geh)$-module generated by the 
extremal weight space $V(\la)_{w\la}$. These modules
$V_w(\la)$ ($w\in W$) are called the Demazure modules.
They are finite dimensional subspaces of $V(\la)$.
Let $(\L(\la),\B(\la))$ be the crystal base of $V(\la)$. In \cite{Ka2}
Kashiwara showed that for each $w\in W$, there exists a subset
$\B_w(\la)$ of $\B(\la)$ such that
\[
\frac{V_w(\la)\cap\L(\la)}{V_w(\la)\cap q\L(\la)}
=\bigoplus_{b\in\B_w(\la)}\Q b.
\]
Furthermore, $\B_w(\la)$ has the following recursive property.
\begin{eqnarray}
&&\mbox{If }r_iw\succ w,\mbox{ then}\nonumber\\
&&\B_{r_iw}(\la)=\bigcup_{n\ge0}\ft{i}^n\B_w(\la)\setminus\{0\}. 
\label{recursive}
\end{eqnarray}
We call $\B_w(\la)$ a {\em Demazure crystal}.

\section{Main Theorem}

Let us present the main theorem in \cite{KMOU} in
the case $\kappa=1$.
For the definition of the {\em mixing index} $\kappa$,
see section 2.3 of \cite{KMOU}.

Let $\la$ be an element of $\Pcll$, and let $B$ be a classical crystal. 
For the theorem, we need to assume four conditions (I-IV).
\begin{description}
\item[(I)   ]
$B$ is perfect of level $l$.
\end{description}
Thus, we can assume an isomorphism between $\B(\la)$ and the set of paths
$\P(\la,B)$. Let $\pbar=\cd\ot\bbar_2\ot\bbar_1$ denote the ground state 
path. Fix a positive integer $d$. For a set of elements
$i_a^{(j)}$ ($j\ge1,1\le a\le d$) in $I$, we define
$B_a^{(j)}$ ($j\ge1,0\le a\le d$) by 
\begin{equation} \label{eq:def_B}
B^{(j)}_0=\{\bbar_j\},\hspace{1cm}
B_a^{(j)}=\bigcup_{n\ge0}\ft{i_a^{(j)}}^n B_{a-1}^{(j)}\setminus\{0\}
\quad(a=1,\cdots,d).
\end{equation}
\begin{description}
\item[(II)  ]
For any $j\ge1$,
$B_d^{(j)}=B$.
\item[(III) ]
For any $j\ge1$ and $1\le a\le d$,
$\langle\la_j,h_{i^{(j)}_a}\rangle\le\veps_{i^{(j)}_a}(b)$
for all $b\in B^{(j)}_{a-1}$.
\end{description}
We now define an element $w^{(k)}$ of the Weyl group $W$ by
\[
w^{(0)}=1,\hspace{1cm}
w^{(k)}=r_{i^{(j)}_a}w^{(k-1)}\quad\mbox{for }k>0,
\]
where $j$ and $a$ are fixed from $k$ by $k=(j-1)d+a,j\ge1,1\le a\le d$.
\begin{description}
\item[(IV)  ]
$w^{(0)}\prec w^{(1)}\prec\cd\prec w^{(k)}\prec\cd$.
\end{description}
See \cite{KMOU}, \cite{KMOTU1} on how to check the last condition.

Now the main statement in \cite{KMOU} is

\begin{theorem}[{\rm \cite{KMOU}}] \label{th:iso}
Under the assumptions (I-IV), we have 
\[
\B_{w^{(k)}}(\la)\simeq u_{\la_j}\ot B^{(j)}_a\ot B^{\ot(j-1)}.
\]
\end{theorem}

The proof is done by showing the recursion relation 
(\ref{recursive}) in the path realization.
The paths on the RHS enjoy 
`full fluctuations' over $B$ in the first $j-1$ steps, while
at the $j$-th step they are allowed only 
`partial fluctuations' over $B^{(j)}_a \subset B$.
After that, they are completely 
frozen to the ground state path $\pbar$.

\section{Example from $\widehat{\mbox{\germlarge sl}}_{\,n}$}

In this section, we describe an example in the $\slnh$ case. We begin with
fixing notations. We use the cyclic notation for $\alpha_i,h_i,\La_i,r_i,
e_i,f_i$, etc, that is, we consider their subscripts $i$ belong to $\Z/n\Z$.
Let $V^{k,l}$ be the irreducible highest weight $U_q(\sln)$-module with 
highest weight $l\La_k$. It turns out that $V^{k,l}$ admits $U'_q(\slnh)$
actions, and has a crystal $B^{k,l}$. 

We describe the explicit actions of $\et{i},\ft{i}$ when $k=1$ (symmetric
tensor case). As a set, $B^{1,l}$ is described as
$B^{1,l}=\{(x_0,\cdots,x_{n-1})\in\Zn^n\mid\sum_{i=0}^{n-1}x_i=l\}$.
The actions of $\et{i},\ft{i}$ are defined as follows.
\begin{eqnarray}
\et{i}(x_0,\cdots,x_{n-1})&=&\left\{
\begin{array}{ll}
(x_0,\cdots,x_{i-1}+1,x_i-1,\cdots,x_{n-1})&(i\ne0)\\
(x_0-1,\cdots,x_{n-1}+1)&(i=0)
\end{array}\right.\label{eq:sln_e}\\
\ft{i}(x_0,\cdots,x_{n-1})&=&\left\{
\begin{array}{ll}
(x_0,\cdots,x_{i-1}-1,x_i+1,\cdots,x_{n-1})&(i\ne0)\\
(x_0+1,\cdots,x_{n-1}-1)&(i=0)
\end{array}\right.\label{eq:sln_f}
\end{eqnarray}
If the right hand side contains a negative component, we should understand
it as 0.

\begin{example}
Let $n=2$. Under the identification of the elements $(2,0)\leftrightarrow00,
(1,1)\leftrightarrow01,(0,2)\leftrightarrow11$, $B^{1,2}$ is described as follows.
\[
B^{1,2}\hspace{5mm}00\mathop{\rightleftharpoons}_0^1
01\mathop{\rightleftharpoons}_0^1 11
\]
Here $b\stackrel{i}{\rightarrow}b'$ means $b'=\ft{i}b$.
\end{example}

We are to show $\la=l\La_0$ and $B=B^{1,l}$ satisfies the four conditions in 
section 3. Firstly, as listed in \cite{KMN2}, $B^{1,l}$ is perfect of level 
$l$. (Note that we deal with the $(A^{(1)}_{n-1},B(l\La_1))$ in their notation.)
$\la_j$ and $\bbar_j$ are given by $\la_j=l\La_{-j}$ and $\bbar_j=
(0,\cdots,0,\stackrel{r}l,0,\cdots,0)$ ($r=-j$ mod $n$,$0\le r\le n-1$). Set
$d=n-1$, and define $i^{(j)}_a\in\Z/n\Z$ by $i^{(j)}_a=a-j$. To illustrate,
we consider the case $j=1$. From the definition (\ref{eq:def_B}) and the rule
(\ref{eq:sln_f}), we easily obtain
\begin{eqnarray*}
&&B^{(1)}_0=\{(0,\cd,0,l)\},
B^{(1)}_1=\{(x_0,0,\cd,0,x_{n-1})\mid x_0+x_{n-1}=l\},\cd,\\
&&B^{(1)}_a=\{(x_0,\cd,x_{a-1},0,\cd,0,x_{n-1})\mid 
x_0+\cd+x_{a-1}+x_{n-1}=l\},\cd.
\end{eqnarray*}
Thus, we have $B^{(1)}_d=B$. For $j>1$, the situation is the same, and $B^{1,l}$
satisfies (II).

\begin{example}
Let $n=3$ and $l=2$. Under the identification of the elements 
$(2,0,0)\leftrightarrow00,(1,1,0)\leftrightarrow01,(1,0,1)\leftrightarrow02,
(0,2,0)\leftrightarrow11,(0,1,1)\leftrightarrow12,(0,0,2)\leftrightarrow22$,
$B^{(1)}_a$ ($a=0,1,2$) are given by
\[
B^{(1)}_0=\{22\},B^{(1)}_1=B^{(1)}_0\sqcup\{00,02\},
B^{(1)}_2=B^{(1)}_1\sqcup\{01,11,12\}(=B).
\]
\end{example}
The third condition is obviously cleared, since 
$\langle\la_j,h_{i^{(j)}_a}\rangle=\langle l\La_{-j},h_{a-j}\rangle=0$ for
$1\le a\le d$. By definition, we have 
\[
w^{(0)}=1,\quad w^{(K)}=\underbrace{r_{K-1}\cd r_1r_0}_K\mbox{ for }K>0,
\]
and the last condition is also true. Therefore, from Theorem \ref{th:iso} we have
\begin{equation} \label{eq:iso1}
\B_{w^{(K)}}(l\La_0)\simeq u_{l\La_{-j}}\ot B^{(j)}_a\ot B^{\ot(j-1)}.
\end{equation}
Here $j$ and $a$ are determined from $K$ by $K=(j-1)d+a,j\ge1,1\le a\le d$,
and $u_{l\La_{-j}}$ is identified with $\cd\ot\bbar_{j+2}\ot\bbar_{j+1}$.

\begin{example} \label{ex:graph}
Let $n=2,l=1$, then $d=1$. Under the identification $(1,0)\leftrightarrow0,
(0,1)\leftrightarrow1$, illustrated in Figure 1 is the Demazure crystal 
$\B_{r_0r_1r_0}(\La_0)$. The symbol $\ot$ is omitted.
\end{example}

\begin{figure}
\centerline{\epsffile{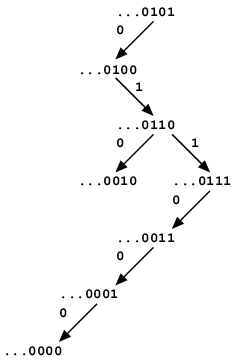}}
\caption{The Demazure crystal $\B_{r_0r_1r_0}(\La_0)$.}
\end{figure}

Let us assume $K=Ld$ for some $L\in\Zn$. Then, (\ref{eq:iso1})
turns out to be
\begin{equation} \label{eq:iso2}
\B_{w^{(Ld)}}(l\La_0)\simeq u_{l\La_{-L}}\ot\left(B^{1,l}\right)^{\ot L}.
\end{equation}
For this fixed $L$, consider a subalgebra $U_q(\sln)=\langle
e_i,f_i,t_i(i\neq-L)\rangle$ of $U_q(\slnh)$. Noting that $\Q(q)u_{l\La_{-L}}$
is a trivial $U_q(\sln)$-module, we see the Demazure crystal 
$\B_{w^{(Ld)}}(l\La_0)$ is isomorphic to $\left(B^{1,l}\right)^{\ot L}$
as $U_q(\sln)$-crystals. This is also true for $q\neq0$, namely, we have
the following theorem.

\begin{theorem} \label{th:sln_inv}
\[
V_{w^{(Ld)}}(l\La_0)\simeq \left(V^{1,l}\right)^{\ot L}
\quad\mbox{as $U_q(\sln)$-modules}.
\]
Here $U_q(\sln)=\langle e_i,f_i,t_i(i\neq-L)\rangle$.
\end{theorem}
As crystals, we already have (\ref{eq:iso2}). Thus it is enough to show
the $U_q(\sln)$-invariance of $V_{w^{(Ld)}}(l\La_0)$. For this purpose,
we cite two propositions from \cite{Ka2}.

\begin{proposition}[{\rm \cite{Ka2} Lemma 3.2.1 (i)}] \label{prop:Ka1}
If $r_iw\prec w$, then 
\[
f_iV_w(\la)\subset V_w(\la).
\]
\end{proposition}

\begin{proposition}[{\rm \cite{Ka2} Corollary 3.2.2}] \label{prop:Ka2}
If $w=r_{i_1}\cd r_{i_l}$ is reduced, then
\[
V_w(\la)=\sum_{k_1,\cd,k_l\ge0}\Q(q)f_{i_1}^{k_1}\cd f_{i_l}^{k_l}u_\la.
\]
\end{proposition}
Let us prove Theorem \ref{th:sln_inv}. We assume $n\ge3$. (The $n=2$ case
is simple.) By definition, $e_iV_{w^{(Ld)}}(l\La_0)\subset 
V_{w^{(Ld)}}(l\La_0)$. We are to show $f_iV_{w^{(Ld)}}(l\La_0)\subset
V_{w^{(Ld)}}(l\La_0)$ for $i\neq-L$. If $r_iw^{(Ld)}\prec w^{(Ld)}$, 
the statement is a direct consequence of Proposition \ref{prop:Ka1}.
(This case includes $i=-L-1$.) Assume $r_iw^{(Ld)}\succ w^{(Ld)}$. 
Since $w^{(Ld)}=r_{Ld-1}\cd r_1r_0$ is reduced, $r_iw^{(Ld)}=
r_ir_{Ld-1}\cd r_1r_0$ should be also reduced. Using the braid relation
$r_jr_{j+1}r_j=r_{j+1}r_jr_{j+1}$ and the condition $i\neq-L,-L-1$,
we can check $r_iw^{(Ld)}=w^{(Ld)}r_{i+L}$. Since $r_iw^{(Ld)}
=r_ir_{Ld-1}\cd r_1r_0$ and $w^{(Ld)}r_{i+L}=r_{Ld-1}\cd r_1r_0r_{i+L}$
are both reduced, from Proposition \ref{prop:Ka2} we have
\begin{eqnarray*}
V_{r_iw^{(Ld)}}(l\La_0)&=&\sum_{k\ge0}f_i^kV_{w^{(Ld)}}(l\La_0),\\
V_{w^{(Ld)}r_{i+L}}(l\La_0)&=&\sum_{k_0,k_1,\cd,k_{Ld}\ge0}\Q(q)
f_{Ld-1}^{k_{Ld}}\cd f_0^{k_1}f_{i+L}^{k_0}u_{l\La_0}\\
&=&\sum_{k_1,\cd,k_{Ld}\ge0}\Q(q)
f_{Ld-1}^{k_{Ld}}\cd f_0^{k_1}u_{l\La_0}\\
&=&V_{w^{(Ld)}}(l\La_0).
\end{eqnarray*}
{}From $V_{r_iw^{(Ld)}}(l\La_0)=V_{w^{(Ld)}r_{i+L}}(l\La_0)$, we can 
conclude the invariance under $f_i$. The theorem is proved.

These facts admit straightforward generalization to arbitrary $k$ cases.
To define the corresponding Weyl group sequence, for $k$ ($1\le k\le n-1$)
and $i\in\Z/n\Z$, we set
\begin{eqnarray*} 
R^{(k)}_i&=&(r_{i+(n-k-1)-(k-1)}\cd r_{i+1-(k-1)}r_{i-(k-1)})\cd\\
&&\qquad\cd(r_{i-1+(n-k-1)}\cd r_ir_{i-1})
(r_{i+(n-k-1)}\cd r_{i+1}r_i).
\end{eqnarray*}
There are $k$ blocks in $R^{(k)}_i$, and in each block there are
$(n-k)$ simple reflections. From the relations among fundamental 
reflections, $R^{(k)}_i$ admits another expression. 
\begin{eqnarray*}
R^{(k)}_i&=&(r_{i-(k-1)+(n-k-1)}\cd r_{i-1+(n-k-1)}r_{i+(n-k-1)})\cd\\
&&\qquad\cd(r_{i+1-(k-1)}\cd r_ir_{i+1})
(r_{i-(k-1)}\cd r_{i-1}r_i).
\end{eqnarray*}
In this case, there are $(n-k)$ blocks, and in each block there are
$k$ simple reflections.
We take $d$ to be $k(n-k)$, and let $w^{(Ld)}$ be determined recursively by
\[
w^{(0)}=1,\quad
w^{((L+1)d)}=R^{(k)}_{-kL}w^{(Ld)}.
\]

\begin{example}
Explicit expression of $w^{(Ld)}$.
\begin{eqnarray*}
n=2,&k=1&w^{(d)}=r_0,w^{(2d)}=r_1r_0,w^{(3d)}=r_0r_1r_0,\cd.\\
n=3,&k=1&w^{(d)}=r_1r_0,w^{(2d)}=r_0r_2r_1r_0,
w^{(3d)}=r_2r_1r_0r_2r_1r_0,\cd.\\
&k=2&w^{(d)}=r_2r_0,w^{(2d)}=r_0r_1r_2r_0,
w^{(3d)}=r_1r_2r_0r_1r_2r_0,\cd.\\
n=4,&k=1&w^{(d)}=r_2r_1r_0,w^{(2d)}=r_1r_0r_3r_2r_1r_0,\\
&&w^{(3d)}=r_0r_3r_2r_1r_0r_3r_2r_1r_0,\cd.\\
&k=2&w^{(d)}=r_0r_3r_1r_0,w^{(2d)}=r_2r_1r_3r_2r_0r_3r_1r_0,\\
&&w^{(3d)}=r_0r_3r_1r_0r_2r_1r_3r_2r_0r_3r_1r_0,\cd.\\
&k=3&w^{(d)}=r_2r_3r_0,w^{(2d)}=r_3r_0r_1r_2r_3r_0,\\
&&w^{(3d)}=r_0r_1r_2r_3r_0r_1r_2r_3r_0,\cd.
\end{eqnarray*}
\end{example}

\begin{theorem} \label{th:rectangular}
We have
\begin{eqnarray*}
\B_{w^{(Ld)}}(l\La_0)&\simeq&u_{l\La_{-kL}}\ot\left(B^{k,l}\right)^{\ot L},\\
V_{w^{(Ld)}}(l\La_0)&\simeq&\left(V^{k,l}\right)^{\ot L}\quad
\mbox{as $U_q(\sln)$-modules},
\end{eqnarray*}
where $U_q(\sln)=\langle e_i,f_i,t_i(i\neq-kL)\rangle$.
\end{theorem}

\section{Demazure Characters and Symmetric Functions}
In this section, we consider the characters of the Demazure modules
we have seen in the previous section. Using the automorphism coming
from the Dynkin diagram symmetry, Theorem \ref{th:rectangular} turns
out to be the following.
\begin{eqnarray}
\B_{w^{(Ld)}}(l\La_{kL})&\simeq&u_{l\La_0}\ot\left(B^{k,l}\right)^{\ot L},
\label{eq:iso3}\\
V_{w^{(Ld)}}(l\La_{kL})&\simeq&\left(V^{k,l}\right)^{\ot L}\quad
\mbox{as $U_q(\sln)$-modules}.\label{eq:iso4}
\end{eqnarray}
Note that $w^{(Ld)}$ is also changed suitably. Here and in what follows,
$U_q(\sln)$ always means the subalgebra of $U_q(\slnh)$ generated by
$e_i,f_i,t_i$ ($i\neq0$).

By definition, the character of the Demazure module $V_w(\la)$ reads as
\[
\ch V_w(\la)=\sum_\mu\sharp\left(\B_w(\la)\right)_\mu e^\mu.
\]
Here $\left(\B_w(\la)\right)_\mu$ is the set of elements in $\B_w(\la)$
of weight $\mu$, and $\mu$ runs over all weights. Consider the character
of $V_{w^{(Ld)}}(l\La_{kL})$ given above. Since it has the 
$U_q(\sln)$-invariance, its character has the following form.
\begin{equation} \label{eq:char}
e^{-l\La_0}\ch V_{w^{(Ld)}}(l\La_{kL})
=\sum_{\la\vdash klL\atop l(\la)\le n}\overline{K}_\la(q)s_\la.
\end{equation}
Here $\la$ runs over all partitions of $klL$ having at most $n$ parts,
$\overline{K}_\la(q)$ is some polynomial in $q$, and $s_\la$ is the 
Schur function considered as a character of $\sln$. $q$ stands for 
$e^{-\delta}$, where $\delta$ is the null root of $\slnh$. In view of
(\ref{eq:iso4}), we have
\[
\left.\left(e^{-l\La_0}\ch V_{w^{(Ld)}}(l\La_{kL})\right)\right|_{q=1}
=s_{(l^k)}^L.
\]
Thus, (\ref{eq:char}) can be viewed as a $q$-analogue of $s_{(l^k)}^L$.

\begin{example}
We consider the case given in Example \ref{ex:graph}.
To adapt the rule in this section, we apply the automorphism of 
the Dynkin diagram. 
\[
e^{-\La_0}\ch V_{r_1r_0r_1}(\La_1)=(1+q)s_{(21)}+q^2s_{(3)}.
\]
\end{example}

Let us examine the polynomial $\overline{K}_\la(q)$. First, we focus
on the case of $k=1$ (symmetric tensor case). The following theorem
was suggested by A.N. Kirillov.

\begin{theorem} \label{th:Kostka}
If $k=1$, we have
\begin{eqnarray*}
\overline{K}_\la(q)&=&q^{-E_0}K_{\la(l^L)}(q),\\
E_0&=&la(L-\frac{n}2(a+1))\quad(a=\left[\frac{L}n\right]).
\end{eqnarray*}
where $K_{\la\mu}(q)$ is the Kostka-Foulkes polynomial.
\end{theorem}
This is a direct consequence of the following expression shown in \cite{NY}.
(See also \cite{DF},\cite{D}.)
\[
K_{\la(l^L)}(q)=\sum q^{\sum_{j=1}^{L-1}jH(b_{j+1}\ot b_j)},
\]
where the sum is over all elements $b_L\ot\cd\ot b_1$ in 
$\left(B^{1,l}\right)^{\ot L}$ which are killed by $\et{i}$ ($i\neq0$) and
have weight $\sum_{i=1}^{n-1}(\la_i-\la_{i+1})\overline{\La}_i$. $H$ stands
for the so-called energy function. Recalling the Milne polynomial
(see p73 of \cite{Ki} and references therein)
\[
M_\mu(x;q)=\sum_\la s_\la(x)K_{\la\mu}(q),
\]
we see the Demazure character (\ref{eq:char}) for $k=1$ turns out to be
the Milne polynomial $M_{(l^L)}(x;q)$ up to a power of $q$.

For general $k$, we have not obtained concrete results yet. We only
mention some generalizations of the Kostka-Foulkes polynomial and
Milne polynomial. For the former, there exists a $q$-analogue of
the multiplicity of the irreducible component $V_\la$ in the
tensor product $V_{\mu_1}\ot\cd\ot V_{\mu_N}$ when all $\mu_i$ have
a rectangular shape. (See (2.35) of \cite{Ki}.) For the latter, there
exists a $q$-analogue of products of Schur functions by Lascoux, Leclerc
and Thibon \cite{LLT}. (They call it H function.) It would be a fascinating
problem to relate them with Demazure characters where $q$ has its own 
meaning as the degree along the null root.

We finish with presenting an inhomogeneous version of (\ref{eq:char})
and Theorem \ref{th:Kostka}.

\begin{theorem} \label{th:inhom}
For a partition $\mu=(\mu_1,\cd,\mu_m)$ of $|\mu|=\mu_1+\cd+\mu_m$,
we set 
\[
w_\mu=R^{(\mu_1)}_{\mu_1}R^{(\mu_2)}_{\mu_1+\mu_2}\cd R^{(\mu_m)}_{|\mu|}.
\]
Then we have 
\[
e^{-\La_0}\ch V_{w_\mu}(\La_{|\mu|})
=q^{c(T_\mu)}\sum_{\la\vdash|\mu|\atop l(\la)\le n}
K_{\la'\mu}(q^{-1})s_\la.
\]
Here $c(T_\mu)$ is the Lascoux-Schuzenberger charge of the tableau
$T_\mu$, whose shape is $(n^pr)$ ($|\mu|=pn+r,0\le r<n$). The contents
of $T_\mu$ $\underbrace{1,1,\cd,1}_{\mu_1},\underbrace{2,2,\cd,2}_{\mu_2},\cd$
are filled in the natural way from left to right from the first row.
$\la'$ stands for the transposition of $\la$.
\end{theorem}
For the definition of the LS charge, see chapter III of \cite{M}.

\begin{example}
Set $n=4$, $\mu=(321)$. Then $|\mu|=6$, $w_\mu=r_1r_2r_3r_1r_0r_2r_1r_0r_3r_2$,
$T_\mu=
\begin{array}{cccc}
1&1&1&2\\
2&3\\
\end{array}$. 
We have
\[
e^{-\La_0}\ch V_{w_\mu}(\La_2)=(1+q)s_{(2^21^2)}
+qs_{(31^3)}+qs_{(2^3)}+q^2s_{(321)}.
\]
\end{example}

We sketch the proof of Theorem \ref{th:inhom}. Firstly, we note 
\[
\B_{w_\mu}(\La_{|\mu|})\simeq 
u_{\La_0}\ot B^{\mu_1,1}\ot B^{\mu_2,1}\ot\cd\ot B^{\mu_m,1}.
\]
This is a consequence of an inhomogeneous version of Theorem \ref{th:iso}.
We also note that $V_{w_\mu}(\La_{|\mu|})$ is $U_q(\sln)$-invariant.
Next, we refer to the following expression of the Kostka-Foulkes polynomial
from \cite{NY}.
\[
K_{\la'\mu}(q^{-1})=\sum_p q^{E(p)},
\]
where $p=u_{\La_0}\ot b_1\ot\cd\ot b_m$ runs over all elements 
in $\B_{w_\mu}(\La_{|\mu|})$ such that $\et{i}p=0$ ($i\neq0$),
$\wt p\equiv\sum_{i=1}^{n-1}(\la_i-\la_{i+1})\overline{\La}_i
\mbox{ mod }\Z\delta$. The energy $E(p)$ of a path 
$p=u_{\La_0}\ot b_1\ot\cd\ot b_m$ is defined by
\[
E(p)=\sum_{j=1}^m\sum_{i=1}^{j-1}
H_{\La_{\mu_i}\La_{\mu_j}}(b_i\ot b_j^{(i+1)}).
\]
Now we have to explain two things: the energy function $H$ and 
the definition of $b_j^{(i+1)}$. Both come from \cite{NY}. Let
us consider the following isomorphism of crystals.
\begin{eqnarray*}
B_1\ot B_2&\simeq&B_2\ot B_1\\
b_1\ot b_2&\mapsto& b'_2\ot b'_1
\end{eqnarray*}
Up to a constant shift, the energy function $H$ on $B_1\ot B_2$ 
is determined by
\begin{eqnarray*}
H(\et{i}(b_1\ot b_2))&=H(b_1\ot b_2)+1
&\mbox{ if }i=0,\et{0}(b_1\ot b_2)=\et{0}b_1\ot b_2,\\
&&\hspace{1.38cm}\et{0}(b'_2\ot b'_1)=\et{0}b'_2\ot b'_1,\\
&=H(b_1\ot b_2)-1
&\mbox{ if }i=0,\et{0}(b_1\ot b_2)=b_1\ot\et{0}b_2,\\
&&\hspace{1.38cm}\et{0}(b'_2\ot b'_1)=b'_2\ot\et{0}b'_1,\\
&\hspace{-.61cm}=H(b_1\ot b_2)
&\mbox{ otherwise. }
\end{eqnarray*}
If $B_1=B^{\mu_i,1}$ and $B_2=B^{\mu_j,1}$, we write 
$H=H_{\La_{\mu_i}\La_{\mu_j}}$. On the other hand,
$b_j^{(i)}$ ($i\le j$) is defined recursively by $b_j^{(j)}=b_j$ and
\begin{eqnarray*}
B^{\mu_i,1}\ot B^{\mu_j,1}&\simeq&B^{\mu_j,1}\ot B^{\mu_i,1},\\
b_i\ot b_j^{(i+1)}&\mapsto&b_j^{(i)}\ot b'_i.
\end{eqnarray*}
Since the crystal graph $\B_{w_\mu}(\La_{|\mu|})$
is connected, the proof reduces to

\begin{proposition}
If $\et{i}p\neq0$, then
\begin{eqnarray*}
E(\et{i}p)&=E(p)-1&(i=0),\\
&\hspace{-.61cm}=E(p)&(i\neq0).
\end{eqnarray*}
\end{proposition}
The case of $i\neq0$ is clear. For $i=0$, we need

\begin{lemma}
Let $p=u_{\La_0}\ot b_1\ot\cd\ot b_m$ and 
$\et{0}p=u_{\La_0}\ot b_1\ot\cd\ot \et{o}b_k\ot\cd\ot b_m\neq0$.
Then we have $k\neq1$, and 
\begin{eqnarray*}
E^{(j)}(\et{0}p)&=E^{(j)}(p)-1&(j=k),\\
&\hspace{-.61cm}=E^{(j)}(p)&(j\neq k).
\end{eqnarray*}
Here $E^{(j)}(p)=\sum_{i=1}^{j-1}H_{\La_{\mu_i}\La_{\mu_j}}
(b_i\ot b_j^{(i+1)})$.
\end{lemma}

\end{document}